%% file: TwoHopsArxiv.tex
\documentclass{IEEEtran}

\usepackage{amsmath,amsfonts,amssymb,euscript, graphicx, 
epsfig,enumerate,float,afterpage, subfigure, ifthen}%
\usepackage[ruled]{algorithm2e}

\newtheorem{theorem}{Theorem}

\newtheorem{lemma}{Lemma}
\newtheorem{definition}{Definition}
 
\newtheorem{assumption}{Assumption} 

\input{SuchaLib.tex}

\begin{document}

\title{Quality of Information Maximization for Wireless Networks via a Fully Separable Quadratic Policy}

\author{
  \IEEEauthorblockN{Sucha Supittayapornpong,~~Michael J. Neely}\\
  \IEEEauthorblockA{Department of Electrical Engineering\\
    University of Southern California\\
    Los Angeles, California\\
    Email: supittay@usc.edu,~~mjneely@usc.edu}
  \thanks{This material was presented in part at the IEEE 
    International Conference on Communications, Ottawa, 
    Canada, June 2012 \cite{Sucha:TwoHop}.

    This material is supported in part by one or more of
    the following: the NSF Career grant CCF-0747525, the
    Network Science Collaborative Technology Alliance sponsored
    by the U.S. Army Research Laboratory W911NF-09-2-0053.}
}

\markboth{arXiv}{}

\maketitle

\begin{abstract}
An information collection problem in a wireless network with random events is considered.  Wireless devices report on each event using one of multiple reporting formats.  Each format has a different quality and uses different data lengths.  Delivering all data in the highest quality format can overload system resources.  The goal is to make intelligent format selection and routing decisions to maximize time-averaged information quality subject to network stability.  Lyapunov optimization theory can be used to solve such a problem by repeatedly minimizing the linear terms of a quadratic drift-plus-penalty expression.  To reduce delays, this paper proposes a novel extension of this technique that preserves the quadratic nature of the drift minimization while maintaining a fully separable structure.  In addition, to avoid high queuing delay, paths are restricted to at most two hops.  The resulting algorithm can push average information quality arbitrarily close to optimum, with a trade-off in queue backlog.  The algorithm compares favorably to the basic drift-plus-penalty scheme in terms of backlog and delay.  Furthermore, the technique is generalized to solve linear programs and yields smoother results than the standard drift-plus-penalty scheme.
\end{abstract}

\section{Introduction}
This paper investigates dynamic scheduling and data format selection in a network where multiple wireless devices, such as smart phones, report information to a receiver station.  The devices together act as a pervasive pool of information about the network environment.  Such scenarios have been recently considered, for example, in applications of social sensing \cite{Miluzzo:Social_sensing} and personal environment monitoring \cite{Kang:SeeMon, Mun:PEIR}.  Sending all information in the highest quality format can quickly overload network resources.  Thus, it is often more important to optimize the \emph{quality of information}, as defined by an end-user, rather than the raw number of bits that are sent.  The case for quality-aware networking is made in \cite{Wang:Beyond_Accuracy, Johnson:QoI,bisdikian-qoi}.  Network management with quality of information awareness for wireless sensor networks is considered in \cite{Liu:QoI_WSN}.  More recently, quality metrics of accuracy and credibility are considered in \cite{BarNoy:QualityOfInformation, Liu:Credibility} using simplified models that do not consider the actual dynamics of a wireless network.  

In this paper, we extend the quality-aware format selection problem in \cite{Liu:Credibility} to a dynamic network setting.  We particularly focus on distributed algorithms for routing, scheduling, and format selection that jointly optimize quality of information.  Specifically, we assume that random events occur over time in the network environment, and these can be sensed by one or more of the wireless devices, perhaps at different sensing qualities.  At the transport layer, each device selects one of multiple reporting formats, such as a video clip at one of several resolution options, an audio clip, or a text message.  Information quality depends on the selected format.  For example, higher quality formats use messages with larger bit lengths.  The resulting bits are handed to the network layer at each device and must be delivered to the receiver station over possibly time-varying channels.  This delivery can be a direct transmission from a device to the receiver station via an uplink channel, or can take a two-hop path that utilizes another device as relay (we restrict paths to at most two-hops for tight control over network delays).  An example is a single-cell wireless network with multiple smart phones and one base station, where each smart phone has 3G capability for uplink transmission and Wi-Fi capability for device-to-device relay transmission.

Such a problem can be cast as a stochastic network optimization and solved using Lyapunov optimization theory.  A ``standard'' method is to minimize a linear term in a quadratic drift-plus-penalty expression \cite{Neely:SNObook, Neely:NOW}.  This can be shown to yield algorithms that converge to optimal average utility with a trade-off in average queue size.  The linearization is useful for enabling decisions to be separated at each device.  However, it can lead to larger queue sizes and delays.  In this work, we propose a novel method that uses a quadratic minimization for the drift-plus-penalty expression, yet still allows separability of the decisions.  This results in an algorithm that maintains distributed decisions across all devices for format selection and routing, similar to the standard (linearized) drift-plus-penalty approach, but reduces overall queue size.

For the derived algorithm, each device observes its input queue length and then selects a format to report an event according to a simple rule.  The routing decision for each group of bits is determined at each device by considering its input, uplink, and relay queues.  Then, allocation of channel resources for direct transmission is determined from a receiver station after observing current uplink queues and channel conditions.  For the relay transmission, an optimization problem involving relay queues, uplink queues and channel conditions is solved at the receiver station to determine an optimal transmission decision.  This process can be decentralized if all channels are orthogonal.

Our analysis shows that the standard drift-plus-penalty algorithm and our new algorithm both converge to the optimal quality of information.  The analysis also shows a deterministic maximum size of each queue.  Simulations show that the new algorithm has a significant savings in queue length which implies reduction of average delay.

Because of the generality of the novel method, it is applied to solve linear programs in the last section.  Linear programs are a special case of the stochastic problems treated in \cite{Neely:NOW}, and hence can be solved by the (linearized) drift plus penalty method of Lyapunov optimization theory.  This is done in \cite{Neely:Dist-Comp} to distributively solve linear programs over graphs.  The current paper applies our novel quadratic drift-plus-penalty algorithm to linear programs to produce smoother results and faster convergence.  Although a solution of this new technique is the time-average of results from multiple iterations, it is different from the ``dual averaging'' method of \cite{Nesterov:DualAveraging} which has a different problem construction, and from the ``alternating direction method of multipliers'' in \cite{Boyd:Alternating} which arises from gradient descent methods rather than from Lyapunov optimization.

Thus, our contributions are threefold:  (i) We formulate an important quality-of-information problem for reporting information in wireless systems.  This problem is of recent interest and can be used in other contexts where ``data deluge'' issues require selectivity in reporting of information.  (ii) We extend Lyapunov optimization theory by presenting a new algorithm that uses a quadratic minimization to reduce queue sizes while maintaining separability across decisions.  This new technique is general and can be used to reduce queue sizes in other Lyapunov optimization problems.  (iii) We illustrate the potential of the quadratic minimization for solving linear programs.

In the next section we formulate the problem.  Sec. \ref{sec:algorithm} derives the novel quadratic algorithm. Sec. \ref{sec:performance} analyzes its performance.  Sec. \ref{sec:simulation} presents simulation results.  Sec. \ref{sec:linear_program} illustrates how to solve linear programs.  The conclusion is in Sec. \ref{sec:conclusion}.

\section{System Model}
\label{sec:model}
\begin{figure}
  \centering
  \includegraphics[scale=0.85]{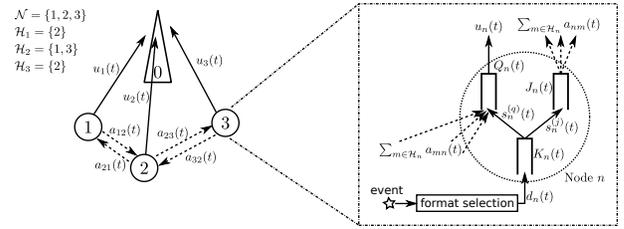}
  \caption{An example network with illustration of the internal queues $K_n(t)$, $Q_n(t)$, $J_n(t)$ for each device $n$.}
  \label{fig:model}
\end{figure}

Consider a network with $N$ wireless devices that report information to a single receiver station.  Let $\set{N} = \{1, \ldots, N\}$ be the set of  devices.  The receiver station is not part of the set $\set{N}$ and can be viewed as ``device 0.'' A network with $N=3$ devices is shown in Fig. \ref{fig:model}.  The system is slotted with fixed size slots $t \in \{0, 1, 2, \ldots\}$.  Every slot, \emph{format selection decisions} are made at the transport layer of each device, and \emph{routing and scheduling} decisions are made at the network layer. 

\subsection{Format Selection}
A new event can occur on each slot.  Events are observed with different levels of quality at each device.  For example, some devices may be physically closer to the event and hence can deliver higher quality.  On slot $t$, each device $n \in \set{N}$ selects a format $f_n(t)$ from a set of available formats $\set{F} = \{0, 1, \ldots, F\}$.  Format selection affects quality and data lengths of the reported information.  To model this, the event on slot $t$ is described by a vector of \emph{event characteristics} $(r_n^{(f)}(t), d_n^{(f)}(t))|_{n \in \set{N}, f \in \set{F}}$.  The value $r_n^{(f)}(t)$ is a numeric \emph{reward} that is earned if device $n$ uses format $f$ to report on the event that occurs on slot $t$.  The value $d_n^{(f)}(t)$ is the amount of data units required for this choice.  This data is injected into the network layer and must eventually be delivered to the receiver station.  To allow a device $n$ not to report on an event, there is a ``blank format'' $0 \in \set{F}$ such that $(r_n^{(0)}(t), d_n^{(0)}(t))=(0,0)$  for all slots $t$ and all devices $n \in \set{N}$.  If a device $n$ does not observe the event on slot $t$ (which might occur if it is physically too far from the event), then $(r_n^{(f)}(t), d_n^{(f)}(t))=(0,0)$ for all formats $f \in \set{F}$.  If no event occurs on slot $t$, then  $(r_n^{(f)}(t), d_n^{(f)}(t)) = (0,0)$ for all $n \in \set{N}$ and $f \in \set{F}$.  

Rewards $r_n(t)$ are assumed to be real numbers that satisfy $0 \leq r_n(t) \leq r_n^\maxvar$ for all $t$, where $r_n^\maxvar$ is a finite maximum.  Data sizes $d_n(t)$ are non-negative integers that satisfy $0 \leq d_n(t) \leq d_n^\maxvar$ for all $t$, where $d_n^\maxvar$ is a finite maximum.  The vectors $(r_n^{(f)}(t), d_n^{(f)}(t))|_{n \in \set{N}, f \in \set{F}}$ are independent and identically distributed (i.i.d.) over slots $t$, and have a joint probability distribution over devices $n$ and formats $f$ that is arbitrary (subject to the above properties). This distribution is not necessarily known.

\subsection{Routing and Scheduling}
\label{ssec:routing_scheduling}
At each device $n \in \set{N}$, the $d_n(t)$ units of data generated by format selection are put into \emph{input queue} $K_n(t)$.  Each device has two orthogonal communication capabilities, called (direct) \emph{uplink transmission} and (ad-hoc) \emph{relay transmission}.  The uplink transmission capability allows each device to communicate to the receiver station directly via an uplink channel.  The relay capability allows communication between a device and its neighboring devices.  To ensure all data takes at most two hops to the destination, the data in each queue $K_n(t)$ is internally routed to one of two queues $Q_n(t)$ and $J_n(t)$, respectively holding data for uplink and relay transmission (see Fig. \ref{fig:model}).  Data in queue $Q_n(t)$ must be transmitted directly to the receiver station, while data in queue $J_n(t)$ can be transmitted to another device $k$, but is then placed in queue $Q_k(t)$ for that device.  This is conceptually similar to the hop-count based queue architecture in \cite{sanjay-shortest-path-ton}. 

In each slot $t$, let $s_n^{(q)}(t)$ and $s_n^{(j)}(t)$ represent the amount of data in $K_n(t)$ that can be internally moved to $Q_n(t)$ and $J_n(t)$, respectively, as illustrated in Fig. \ref{fig:model}.  These decision variables are chosen within sets $\set{S}_n^{(q)}$ and $\set{S}_n^{(j)}$, respectively, where: 
\begin{eqnarray*}
  \set{S}_n^{(q)} &\defequiv& \{0, 1, \dotsc, s_n^{(q)\maxvar}\} \\
  \set{S}_n^{(j)} &\defequiv& \{0, 1, \dotsc, s_n^{(j)\maxvar}\}
\end{eqnarray*}
where $s_n^{(q)\maxvar}$, $s_n^{(j)\maxvar}$ are finite maximum values. 
Then the dynamics of $K_n(t)$ are: 
\begin{equation} 
  \label{eq:input_queue}
  K_n(t+1) = \max[K_n(t) - s_n^{(q)}(t) - s_n^{(j)}(t), 0] + d_n(t) 
\end{equation} 
As a minor technical detail that is useful later, the $\max[\dotsb, 0]$ operation above allows the $s_n^{(q)}(t)$ and $s_n^{(j)}(t)$ decisions to sum to more than $K_n(t)$.  The \emph{actual} $s_n^{(q)\actvar}(t)$ and $s_n^{(j)\actvar}(t)$ data units moved from $K_n(t)$ can be any values that satisfy: 
\begin{gather} 
  s_n^{(q)\actvar}(t) + s_n^{(j)\actvar}(t) = \min[K_n(t), s_n^{(q)}(t) + s_n^{(j)}(t)] \label{eq:actual_sq_sj}\\
  0 \leq s_n^{(q)\actvar}(t) \leq s_n^{(q)}(t) \label{eq:actual_sq}\\
  0 \leq s_n^{(j)\actvar}(t) \leq s_n^{(j)}(t) \label{eq:actual_sj}
\end{gather} 

Wireless transmission is assumed to be channel-aware, and decision options are determined by a vector $\vect{\eta}(t)$ of \emph{current channel states} in the network.  Specifically, let $u_n(t)$ be the amount of uplink data that can be transmitted from device $n$ to the receiver station, and let $\vect{u}(t) = (u_n(t))|_{n \in \set{N}}$ be the vector of these transmission decisions.  It is assumed that $\vect{u}(t)$ is chosen every slot $t$ within a set $\set{U}_{\vect{\eta}(t)}$ that depends on the observed $\vect{\eta}(t)$.  Similarly, let $a_{nm}(t)$ be the amount of data selected for ad-hoc transmission between devices $n$ and $m$, and let $\vect{a}(t)  = (a_{nm}(t))|_{n,m \in \set{N}}$ and $a_{nn}(t) = 0$ for every $t$ and $n$.  These transmissions are assumed to be orthogonal to the uplink transmissions.  Every slot $t$, the $\vect{a}(t)$ vector is chosen within a set $\set{A}_{\vect{\eta}(t)}$ that depends on the observed $\vect{\eta}(t)$.  The sets $\set{U}_{\vect{\eta}(t)}$ and $\set{A}_{\vect{\eta}(t)}$ depend on the resource allocation, modulation, and coding options for transmission.  If each uplink channel is orthogonal then set $\set{U}_{\vect{\eta}(t)}$ can be decomposed into a set product of individual options for each uplink, where each option depends on the component of $\vect{\eta}(t)$ that represents its own uplink channel.  Orthogonal relay links can be treated similarly. 

The dynamics of relay queue $J_n(t)$ are: 
\begin{equation}
  \label{eq:relay_queue}
  J_n(t+1) = \max \prts{ J_n(t) - \mbox{$\sum_{m \in \set{N}}$} a_{nm}(t) + s_n^{(j)\actvar}(t) , 0 }.
\end{equation}
As before, the \emph{actual} amount of data $a_{nm}^\actvar(t)$ satisfies:
\begin{gather}
  \mbox{$\sum_{m \in \set{N}}$} a_{nm}^\actvar(t) = \min \Bigl( J_n(t) + s_n^{(j)\actvar}(t), \mbox{$\sum_{m \in \set{N}}$} a_{nm}(t) \Bigr) \label{eq:actual_transmit} \\
  0 \leq a_{nm}^\actvar(t) \leq a_{nm}(t) \quad\quad \text{for}~ m \in \set{N} \label{eq:actual_relay}.
\end{gather}

The dynamics of uplink queue $Q_n(t)$ are:
\begin{multline}
  \label{eq:uplink_queue}
  Q_n(t+1) = \max \prts{ Q_n(t) - u_n(t) + s_n^{(q)\actvar}(t) , 0 }\\ 
  + \mbox{$\sum_{m \in \set{N}}$} a_{mn}^\actvar(t).
\end{multline}
Notice that all data transmitted to a relay is placed in the uplink queue of that relay (which ensures all paths take at most two hops).  The queueing equations \eqref{eq:relay_queue} and \eqref{eq:uplink_queue} involve actual amounts of data, but they can be bounded using \eqref{eq:actual_sq}, \eqref{eq:actual_sj} and \eqref{eq:actual_relay} as
\begin{align}
  J_n(t+1) & \leq \max \prts{ J_n(t) - \mbox{$\sum_{m \in \set{N}}$} a_{nm}(t) + s_n^{(j)}(t) , 0 }  \label{eq:relay_queue_bound} \\
  Q_n(t+1) & \leq \max \prts{ Q_n(t) - u_n(t) + s_n^{(q)}(t) , 0 } \notag \\
  & \hspace{10em} + \mbox{$\sum_{m \in \set{N}}$} a_{mn}(t).  \label{eq:uplink_queue_bound}
\end{align}

The queue dynamics \eqref{eq:input_queue}, \eqref{eq:relay_queue_bound}, \eqref{eq:uplink_queue_bound} do not require the actual variables $s_n^{(j)\actvar}, s_n^{(q)\actvar}(t)$, $a_{nm}^\actvar(t)$, and are the only ones needed in the rest of the paper.

Assume the decision sets $\set{U}_{\vect{\eta}(t)}$ and $\set{A}_{\vect{\eta}(t)}$ ensure that transmissions have bounded rates.  Specifically, let $u_n^\maxvar$ and $a_{nm}^\maxvar$ be finite maximum values of $u_n(t)$ and $a_{nm}(t)$.  Further, assume that for each $n \in \set{N}$, $s_n^{(q)\maxvar} \geq u_n^{\maxvar}$ and $s_n^{(j)\maxvar} \geq \sum_{m\in\set{N}}a_{nm}^\maxvar$, so that the maximum amount that can be internally shifted is at least as much as the maximum amount that can be transmitted.


\subsection{Stochastic Network Optimization}
\label{ssec:stochastic_optimization}
Here we define the problem of maximizing time-averaged quality of information subject to queue stability.  We use the following definitions \cite{Neely:NOW}: 
\begin{definition}
  \label{def:queue_stability}
  Queue $\prtc{ X(t) : t \in \prtc{0, 1, 2, \dotsc}}$ is strongly stable if
  \begin{equation*}
    \limsup_{t \rightarrow \infty} \mbox{$\frac{1}{t} \sum_{\tau = 0}^{t-1}$} \expect{ X(\tau) } < \infty
  \end{equation*}
\end{definition}

\begin{definition}
  \label{def:network_stability}
  A network of queues is strongly stable if every queue in the network is strongly stable.
\end{definition}

In words, definition \ref{def:queue_stability} means that a queue is strongly stable if its average queue backlog is finite.

Let $y_0(t) \defequiv \sum_{n \in \set{N}} r_n(t)$ be the total quality of information from format selection on slot $t$, and $y_0^\maxvar \triangleq \sum_{n \in \set{N}} r_n^\maxvar$ is its upper bound.  The time-averaged total information quality is
\begin{equation*}
  \bar{y}_0 \defequiv \liminf_{t \rightarrow \infty} \mbox{$\frac{1}{t} \sum_{\tau = 0}^{t-1}$} \expect{y_0(\tau)}.
\end{equation*}

For simplicity of notation, let $\random(t)$ represent a collective vector of event and channel randomness on slot $t$, and let $\control(t)$ be a collective vector of all decision variables on slot $t$: 
\begin{eqnarray*}
  \random(t) &\defequiv& [\vect{\eta}(t); (r_n^{(f)}(t), d_n^{(f)}(t))|_{n\in\set{N}, f\in\set{F}}] \\
  \control(t) &\defequiv& [\vect{a}(t); \vect{u}(t); (f_n(t))|_{n\in\set{N}}; (s_n^{(q)}(t), s_n^{(j)}(t))|_{n \in \set{N}}] 
\end{eqnarray*}

It is our objective to solve: 
\begin{align}
  \label{eq:optimization_problem}
  \maximize \quad
  & \bar{y}_0 \\
  \subjectto \quad
  & \text{Network is strongly stable} \notag \\
  & \control(t) \in \Phi_{\random(t)} ~ \text{for all} ~ t, \notag
\end{align}
where $\Phi_{\random(t)}$ is a feasible set of control actions depending on randomness at time $t$. So, any selected $\control(t) \in \Phi_{\random(t)}$ yields:
\begin{align*}
  & f_n(t) \in \set{F} ~ \text{for all} ~ n \in \set{N} \\
  & s_n^{(q)}(t) \in \set{S}_n^{(q)} ~ \text{for all} ~ n \in \set{N} \\
  & s_n^{(j)}(t) \in \set{S}_n^{(j)} ~ \text{for all} ~ n \in \set{N} \\
  & \vect{u}(t) \in \set{U}_{\vect{\eta}(t)} \\
  & \vect{a}(t) \in \set{A}_{\vect{\eta}(t)} 
\end{align*}
This problem is always feasible because stability is trivially achieved if all devices always select the blank format.

\section{Dynamic Algorithm}
\label{sec:algorithm}
This section derives a novel ``quadratic policy'' to solve problem \eqref{eq:optimization_problem}.  
The policy gives faster convergence and smaller queue sizes as compared to the 
``standard'' drift-plus-penalty (or ``max-weight'') policy of \cite{Neely:SNObook, Neely:NOW}.

\subsection{Lyapunov Optimization}
Let $\queue(t) = ( K_n(t), Q_n(t), J_n(t) )|_{n \in \set{N}}$ represent a vector of all queues in the system.

Define a quadratic \emph{Lyapunov function} $L(\queue(t)) \defequiv \frac{1}{2} \sum_{n \in \set{N}} \left[ K_n^2(t) + Q_n^2(t) + J_n^2(t) \right]$.  Then the Lyapunov drift, the difference of Lyapunov functions between two consecutive slots, is defined by $L(\queue(t+1)) - L(\queue(t))$.  

In order to maximize $\bar{y}_0$ in \eqref{eq:optimization_problem}, the drift-plus-penalty function $L(\queue(t+1)) - L(\queue(t)) - V y_0(t)$ is considered, where $V \geq 0$ is a constant that determines a trade-off between queue size and proximity to the optimality.\footnote{The minus sign in front of $V y_0(t)$ is because the quality of information can be viewed as a negative penalty.}

Later, this drift is used to show stability of a system.  Intuitively, when queue lengths grow large beyond certain values, the drift becomes negative and a system is stable because the negative drift roughly implies reduction of total queue lengths.

Let $\RealSet$ and $\RealSet_+$ denote the set of real numbers and non-negative real numbers, respectively.

\begin{lemma}
  \label{lem:absbounds}
Let $a_i \in \RealSet$ and $b_j \in \RealSet_+$ for $i \in \{0, 1, 2, \dotsc, A \}$ and $j \in \{0, 1, 2, \dotsc, B \}$.  Assume further that $\abs{a_i} \leq a_i^\maxvar$ and $\abs{b_j} \leq b_j^\maxvar$ for each feasible $i$ and $j$.  Then for any $x \in \RealSet_+$,
\begin{align}
  & \prts{ \max\prtr{ x + \mbox{$\sum_{i = 1}^A$} a_i , 0 } + \mbox{$\sum_{j = 1}^B$} b_j }^2 - x^2 \notag \\
  & \leq \mbox{$\sum_{i = 1}^A$} ( x + a_i )^2 + \mbox{$\sum_{j = 1}^B$} ( x + b_j )^2 - (A + B )x^2 + C     \label{eq:quadratic_bound} \\
  & \leq 2x\prtbs{ \mbox{$\sum_{i = 1}^A$} a_i + \mbox{$\sum_{j = 1}^B$} b_j } + C'  \label{eq:linear_bound}
\end{align}
where
\begin{align*}
  C & = 2 \left[ \mbox{$\sum_{i = 1}^{A}\sum_{i'=1}^{i-1}$} a_i^\maxvar a_{i'}^\maxvar + \mbox{$\sum_{j = 1}^B\sum_{j' = 1}^{j - 1}$} b_j^\maxvar b_{j'}^\maxvar \right. \\
  & \hspace{12.5em} \left.  + \mbox{$\sum_{i = 1}^A \sum_{j = 1}^B$} a_i^\maxvar b_j^\maxvar \right] \\
  C' & = \prts{ \mbox{$\sum_{i = 1}^A$} a_i^\maxvar + \mbox{$\sum_{j = 1}^B$} b_j^\maxvar }^2
\end{align*}
Note that the first bound \eqref{eq:quadratic_bound} is used in the quadratic policy, while the second bound \eqref{eq:linear_bound} can lead to the max-weight policy.
\end{lemma}

\begin{IEEEproof}
  \begin{align}
    & \prts{ \max\prtr{ x + \mbox{$\sum_{i = 1}^A$} a_i , 0 } + \mbox{$\sum_{j = 1}^B$} b_j }^2 - x^2 \notag \\
    & {\scriptstyle \leq \prtr{ x + \sum_{i = 1}^A a_i }^2 + \prtr{ \sum_{j = 1}^B b_j }^2 + 2 \sum_{j = 1}^B b_j \abs{ x + \sum_{i = 1}^A a_i } - x^2 }\notag \\
    & {\scriptstyle = 2x \sum_{i = 1}^A a_i + \prtr{ \sum_{i = 1}^A a_i }^2 + \prtr{ \sum_{j = 1}^B b_j }^2 + 2 \sum_{j = 1}^B b_j \abs{ x + \sum_{i = 1}^A a_i } }\notag \\
    & {\scriptstyle \leq 2x \sum_{i = 1}^A a_i + \sum_{i = 1}^A a_i^2 + 2 \sum_{i = 1}^A \sum_{i' = 1}^{i - 1} \abs{a_i a_{i'}} + \sum_{j = 1}^B b_j^2 } \notag \\
    & \quad {\scriptstyle + 2 \sum_{j = 1}^B\sum_{j' = 1}^{j - 1} b_j b_{j'} + 2\sum_{j = 1}^B b_j \abs{ x + \sum_{i = 1}^A \abs{a_i} } } \notag \\
    & {\scriptstyle = 2x \sum_{i = 1}^A a_i + \sum_{i = 1}^A a_i^2 + 2\sum_{i = 1}^A\sum_{i' = 1}^{i - 1} \abs{a_i a_{i'}} + \sum_{j = 1}^B b_j^2 } \notag \\
    & \quad {\scriptstyle + 2\sum_{j = 1}^B\sum_{j' = 1}^{j - 1} b_j b_{j'} + 2\sum_{j = 1}^B b_j x + 2\sum_{j = 1}^B \sum_{i = 1}^A b_j\abs{a_i} } \notag \\
    & {\scriptstyle = \sum_{i = 1}^A \prtr{ x + a_i }^2 + \sum_{j = 1}^B \prtr{x + b_j}^2 - (A + B)x^2 } \notag \\
    & \quad {\scriptstyle + 2\sum_{i = 1}^A\sum_{i' = 1}^{i - 1} \abs{a_i a_{i'}} + 2\sum_{j = 1}^B\sum_{j' = 1}^{j - 1} b_j b_{j'} + 2\sum_{i = 1}^A \sum_{j = 1}^B \abs{a_i}b_j }\notag \\
    & {\scriptstyle \leq \sum_{i = 1}^A \prtr{ x + a_i }^2 + \sum_{j = 1}^B \prtr{x + b_j}^2 - (A + B)x^2 + C } \label{eq:bound_max_bp} \\ 
    & {\scriptstyle \leq 2x \prts{ \sum_{i = 1}^A a_i + \sum_{j = 1}^B b_j } + \sum_{i = 1}^A a_i^{\maxvar 2} + \sum_{j = 1}^B b_j^{\maxvar 2} + C } \notag \\
    & {\scriptstyle = 2x \prts{ \sum_{i = 1}^A a_i + \sum_{j = 1}^B b_j } + C' } \label{eq:bound_max_bp_linear}
  \end{align}
Inequalities \eqref{eq:bound_max_bp} and \eqref{eq:bound_max_bp_linear} prove respectively relation \eqref{eq:quadratic_bound} and \eqref{eq:linear_bound}.
\end{IEEEproof}

Using queuing dynamic \eqref{eq:input_queue}, \eqref{eq:relay_queue_bound}, and \eqref{eq:uplink_queue_bound}, the drift-plus-penalty is bounded by \eqref{eq:pure_drift} below.  Then, using relation \eqref{eq:quadratic_bound}, the bound becomes \eqref{eq:square_drift}. 
\begin{align}
  & L(\queue(\tau+1)) - L(\queue(\tau)) - V y_0(\tau) \notag \\
  & \scriptstyle \leq \frac{1}{2} \sum_{n \in \set{N}} \bigl\{ \left[ \max( K_n(\tau) - s_n^{(q)}(\tau) - s_n^{(j)}(\tau), 0 ) + d_n(\tau) \right]^2 - K_n(\tau)^2 \notag \\
  & \scriptstyle \quad + \left[ \max( Q_n(\tau) - u_n(\tau) + s_n^{(q)}(\tau), 0 ) + \sum_{m \in \set{N}} a_{mn}(\tau) \right]^2 - Q_n(\tau)^2 \notag \\
  & \scriptstyle \quad + \left[ \max( J_n(\tau) - \sum_{m \in \set{N}} a_{nm}(\tau) + s_n^{(j)}(\tau), 0 ) \right]^2 - J_n(\tau)^2 - 2V r_n(\tau) \bigr\} \label{eq:pure_drift} \\
  & \scriptstyle{\leq \frac{1}{2} \sum_{n \in \set{N}} \bigl\{ \left[ K_n(\tau) - s_n^{(q)}(\tau) \right]^2 + \left[ K_n(\tau) - s_n^{(j)}(\tau) \right]^2 + \left[ K_n(\tau) + d_n(\tau) \right]^2 } \notag \\
  & \scriptstyle{\quad + \left[Q_n(\tau) - u_n(\tau) \right]^2 + \left[ Q_n(\tau) + s_n^{(q)}(\tau) \right]^2 + \sum_{m \in \set{N}} \left[Q_n(\tau) + a_{mn}(\tau) \right]^2 } \notag \\
  & \scriptstyle{\quad + \sum_{m \in \set{N}} \left[J_n(\tau) - a_{nm}(\tau) \right]^2 + [ J_n(\tau) + s_n^{(j)}(\tau) ]^2 - 2V r_n(\tau) + D_n(\tau) \bigr\} } \label{eq:square_drift}
\end{align}
where
\begin{align*}
\scriptstyle D_n(\tau) \triangleq 
& \scriptstyle -3 K_n^2(\tau) - (2 + \abs{\set{N}}) Q_n^2(\tau) - (1 + \abs{\set{N}}) J_n^2(\tau) \\
& \scriptstyle + 2 s_n^{(q)\maxvar} s_n^{(j)\maxvar} + 2 s_n^{(q)\maxvar} d_n^\maxvar + 2 s_n^{(j)\maxvar} d_n^\maxvar \\
& \scriptstyle + 2 u_n^\maxvar s_n^{(q)\maxvar} + 2 u_n^\maxvar \sum_{m \in \set{N}} a_{mn}^\maxvar + 2 s_n^{(q)\maxvar} \sum_{m \in \set{N}} a_{mn}^\maxvar \\
& \scriptstyle + \sum_{m \in \set{N}} \sum_{m' \in \set{N} - \{m\}} a_{mn}^\maxvar a_{m'n}^\maxvar + 2 s_n^{(j)\maxvar} \sum_{m \in \set{N}} a_{nm}^\maxvar \\
& \scriptstyle + \sum_{m \in \set{N}} \sum_{m' \in \set{N} - \{m\}} a_{nm}^\maxvar a_{nm'}^\maxvar
\end{align*}

Minimizing the actual drift-plus-penalty term \eqref{eq:pure_drift} is computationally expensive.  In this paper, we propose a novel \textit{quadratic policy}, derived from \eqref{eq:square_drift}, that preserves the quadratic nature of the actual minimization while keeping decisions separable.  As a result, the policy leads to a separated control algorithm in Sec. \ref{ssec:separability}.

\begin{definition}
  \label{def:quadratic_policy}
  Every time $t$, the quadratic policy observes current queue backlogs $\queue(t)$ and randomness $\random(t)$. Then it makes a decision according to the following minimization problem.
  \begin{align*}
    \minimize \quad
    & \scriptstyle \sum_{n \in \set{N}} \bigl\{ \left[ K_n(t) - s_n^{(q)}(t) \right]^2 + \left[ K_n(t) - s_n^{(j)}(t) \right]^2  \notag \\
    & \scriptstyle + \left[ K_n(t) + d_n(t) \right]^2 + \left[Q_n(t) - u_n(t) \right]^2 + \left[ Q_n(t) + s_n^{(q)}(t) \right]^2 \notag \\
    & \scriptstyle + \sum_{m \in \set{N}} \left[Q_n(t) + a_{mn}(t) \right]^2 + \sum_{m \in \set{N}} \left[J_n(t) - a_{nm}(t) \right]^2 \notag \\
    & \scriptstyle + \left[ J_n(t) + s_n^{(j)}(t) \right]^2 - 2V r_n(t) \bigr\} \\
    \subjectto \quad
    & \scriptstyle s_n^{(q)}(t) \in \set{S}_n^{(q)}, ~ s_n^{(j)}(t) \in \set{S}_n^{(j)} ~~~ \forall n \in \set{N} \\
    & \scriptstyle f_n(t) \in \set{F}, d_n(t) = d_n^{(f_n(t))}(t), ~ r_n(t) = r_n^{(f_n(t))}(t)  ~~~ \forall n \in \set{N} \\
    & \scriptstyle \vect{a}(t) \in \set{A}_{\vect{\eta}(t)}, ~ \vect{u}(t) \in \set{U}_{\vect{\eta}(t)} 
  \end{align*}
\end{definition}

\subsection{Separability}
\label{ssec:separability}
The control algorithm can be derived from the quadratic policy in definition \ref{def:quadratic_policy}.  The whole minimization can be done separately due to a unique structure of the quadratic policy.  This leads to five subproblems, as described below.

At every slot $t$, each device $n \in \set{N}$ observes input queue $K_n(t)$ and options $(r_n^{(f)}(t), d_n^{(f)}(t))|_{f \in \set{F}}$. It then chooses a format $f_n(t)$ according to the \emph{admission-control problem}:
\begin{align}
  \label{eq:admission_control}
  \minimize \quad
  & \prts{K_n(t) + d_n^{(f_n(t))}(t)}^2 - 2V r_n^{(f_n(t))}(t) \\
  \subjectto \quad
  & f_n(t) \in \set{F} \notag
\end{align}
This is solved easily by comparing each option $f_n(t) \in \set{F}$. 

Each device $n$ moves data from its input queue to its uplink queue according to the \emph{uplink routing problem}
\begin{align}
  \label{eq:uplink_buffer_control}
  \minimize \quad
  &  \prts{K_n(t) - s_n^{(q)}(t)}^2 + \prts{Q_n(t) + s_n^{(q)}(t)}^2 \\
  \subjectto \quad
  & s_n^{(q)}(t) \in \set{S}_n^{(q)} \notag.
\end{align}
This can be solved in a closed form by letting $I_{Q}^{+}(t) \triangleq \bigl\lceil \frac{K_n(t) - Q_n(t)}{2} \bigr\rceil$, $I_{Q}^{-}(t) \triangleq \bigl\lfloor \frac{K_n(t) - Q_n(t)}{2} \bigr\rfloor$ and $g_Q(x, t) = \prts{K_n(t) - x}^2 + \prts{Q_n(t) + x}^2$.  Then choose
\begin{align}
  \label{eq:closed_form_snq}
  & s_n^{(q)}(t) = \\
  & \left\{
  \begin{array}{ll}
    \scriptstyle s_n^{(q)\maxvar}  &,~ \scriptstyle K_n(t) - Q_n(t) \geq 2 s_n^{(q)\maxvar} \\
    \scriptstyle \argmin_{x \in \left\{I_{Q}^{+}(t), I_{Q}^{-}(t) \right\} } g_Q(x,t)  &,~ \scriptstyle 0 < K_n(t) - Q_n(t) < 2 s_n^{(q)\maxvar} \\
    \scriptstyle 0  &,~ \scriptstyle K_n(t) - Q_n(t) \leq 0
  \end{array}
  \right. \notag
\end{align}

Also each device $n$ moves data from its input queue to its relay queue according to the \emph{relay routing problem}
\begin{align}
  \label{eq:relay_buffer_control}
  \minimize \quad
  &  \prts{K_n(t) - s_n^{(j)}(t)}^2 + \prts{J_n(t) + s_n^{(j)}(t)}^2. \\
  \subjectto \quad
  & s_n^{(j)}(t) \in \set{S}_n^{(j)} \notag
\end{align}
Again, let $I_{J}^{+}(t) \triangleq \bigl\lceil \frac{K_n(t) - J_n(t)}{2} \bigr\rceil$, $I_{J}^{-}(t) \triangleq \bigl\lfloor \frac{K_n(t) - J_n(t)}{2} \bigr\rfloor$ and $g_J(x, t) = \prts{K_n(t) - x}^2 + \prts{J_n(t) + x}^2$.  Then choose
\begin{align}
  \label{eq:closed_form_snj}
  & s_n^{(j)}(t) = \\
  & \left\{
  \begin{array}{ll}
    \scriptstyle s_n^{(j)\maxvar}  &,~ \scriptstyle K_n(t) - J_n(t) \geq 2 s_n^{(j)\maxvar} \\
    \scriptstyle \arg\min_{x \in \{ I_{J}^+(t), I_{J}^-(t) \}} g_J(x, t)  &,~ \scriptstyle 0 < K_n(t) - J_n(t) < 2 s_n^{(j)\maxvar} \\
    \scriptstyle 0  &,~ \scriptstyle K_n(t) - J_n(t) \leq 0
  \end{array}
  \right. \notag
\end{align}
Note that the solutions from the quadratic policy are ``smoother'' as compared to the solutions from the max-weight policy that would choose ``bang-bang'' decisions of either $0$ or $s_n^{(q)\maxvar}$ for $s_n^{(q)}(t)$ (and $0$ or $s_n^{(j)\maxvar}$ for $s_n^{(j)}(t)$).

The \emph{uplink allocation} problem is
\begin{align}
  \label{eq:uplink_resource}
  \minimize \quad
  & \mbox{$\sum_{n \in \set{N}}$} \prts{ Q_n(t) - u_n(t) }^2 \\
  \subjectto \quad
  & \vect{u}(t) \in \set{U}_{\vect{\eta}(t)}. \notag
\end{align}
This can be solved at the receiver station. If all uplink channels are orthogonal, the problem can be decomposed further to be solved at each device $n$ by
\begin{align}
  \label{eq:uplink_resource_orth}
  \minimize \quad
  &  \prts{Q_n(t) - u_n(t)}^2 \\
  \subjectto \quad
  & u_n(t) \in \set{U}_{n,\vect{\eta}(t)}, \notag
\end{align}
where $\set{U}_{n,\vect{\eta}(t)}$ is a feasible set of $u_n(t)$.  An optimal uplink transmission rate is the closest rate in $\set{U}_{n, \vect{\eta}(t)}$ to $Q_n(t)$.

The \emph{relay allocation} problem is
\begin{align}
  \label{eq:relay_resource}
  \minimize \quad
  &  \mbox{$\sum_{n \in \set{N}} \sum_{m \in \set{N}}$} \Bigl\{ \prts{Q_n(t) + a_{mn}(t)}^2 \Bigr. \notag\\
  &  \hspace{7em} \Bigl. + \prts{J_n(t) -  a_{nm}(t)}^2 \Bigr\} \\
  \subjectto \quad
  & \vect{a}(t) \in \set{A}_{\vect{\eta}(t)}. \notag
\end{align}
If channels are orthogonal so the sets have a product form, then the decisions are separable across transmission links $(n,m)$ for $n \in \set{N}, m \in \set{N}$ as
\begin{align}
  \label{eq:relay_resource_orth}
  \minimize \quad 
  & \prts{ Q_m(t) + a_{nm}(t) }^2 + \prts{ J_n(t) -  a_{nm}(t) }^2 \\
  \subjectto \quad
  & a_{nm}(t) \in \set{A}_{nm,\vect{\eta}(t)}, \notag
\end{align}
where $\set{A}_{nm,\vect{\eta}(t)}$ is a feasible set of $a_{nm}(t)$.  The closed form solution of this problem is
\begin{align}
  \label{eq:closed_form_anmt}
  & a_{nm}(t) = \\
  & \left\{
  \begin{array}{ll}
    \scriptstyle a_{nm}^\maxvar 	&, \scriptstyle ~ J_n(t) - Q_m(t) \geq 2 a_{nm}^\maxvar \\
    \scriptstyle \argmin_{x \in \{ I_{A}^+(t), I_{A}^-(t) \}} g_A(x, t) &, \scriptstyle ~ 0 < J_n(t) - Q_m(t) < 2 a_{nm}^\maxvar \\
    \scriptstyle 0					&, \scriptstyle ~ J_n(t) - Q_m(t) \leq 0
  \end{array}
  \right. \notag
\end{align}
where $I_{A}^{+}(t) \triangleq \argmin_{a \in \set{A}_{nm},\vect{\eta}(t)} \abs{ a - \frac{J_n(t) - Q_m(t)}{2} }$ and $I_{A}^{-}(t) \triangleq \argmin_{a \in \set{A}_{nm,\vect{\eta}(t)}-\{I_{A}^{+}(t)\}} \abs{ a - \frac{J_n(t) - Q_m(t)}{2} }$ and $g_A(x, t) = \prts{J_n(t) - x}^2 + \prts{Q_m(t) + x}^2$.

\subsection{Algorithm}
\label{ssec:algorithm}
At every time slot $t$, our algorithm has two parts: device side and receiver-station side.
\begin{algorithm}
  \hspace{-1em}\tcp{Device side}
  \ForEach{device $n \in \set{N}$}{
    -- Observe $K_n(t), Q_n(t)$ and $J_n(t)$ \\
    -- Observe $(r_n^{(f)}(t), d_n^{(f)}(t))|_{f \in \set{F}}$ \\
    -- Select format $f_n(t)$ according to \eqref{eq:admission_control} \\
    -- Move data from $K_n(t)$ to $Q_n(t)$ and $J_n(t)$ with $s_n^{(q)\actvar}(t), s_n^{(j)\actvar}(t)$ satisfying \eqref{eq:actual_sq_sj}-\eqref{eq:actual_sj} and \eqref{eq:actual_transmit}-\eqref{eq:actual_relay} with values of $s_n^{(q)}(t), s_n^{(j)}(t)$ calculated from \eqref{eq:closed_form_snq} and \eqref{eq:closed_form_snj}.
  }
  \caption{Distributed format selection and routing}
\end{algorithm}

\begin{algorithm}
  \hspace{-1em}\tcp{Receiver-station side}
  \For{receiver station $0$}{
    -- Observe $(Q_n(t), J_n(t))|_{n \in \set{N}}$ \\
    -- Observe $\set{U}_{\vect{\eta}(t)}$ and $\set{A}_{\vect{\eta}(t)}$ \\
    -- Signal devices $n \in \set{N}$ to make uplink transmission $\vect{u}(t)$ according to \eqref{eq:uplink_resource} \\
    -- Signal devices $n \in \set{N}$ to relay data $\vect{a}(t)$ according to \eqref{eq:relay_resource} 
  }
  \caption{Uplink and Relay resource allocation}
\end{algorithm}

%
After these processes, queues $K_n(t+1), Q_n(t+1)$ and $J_n(t+1)$ are updated via \eqref{eq:input_queue}, \eqref{eq:relay_queue}, \eqref{eq:uplink_queue}.

\section{Stability and Performance Bounds}
\label{sec:performance}
Compare the quadratic policy with any other policy.  Let $(f_n(\tau), s_n^{(q)}(\tau), s_n^{(j)}(\tau))|_{n \in \set{N}}, \vect{u}(\tau), \vect{a}(\tau)$ be the decision variables from the quadratic policy in definition \ref{def:quadratic_policy}. From $f_n(\tau)$, $r_n(t) \triangleq r_n^{(f_n(t))}(t)$ and $d_n(t) \triangleq d_n^{(f_n(t))}(t)$.  Then, let $(\hat{f}_n(\tau), \hat{s}_n^{(q)}(\tau), \hat{s}_n^{(j)}(\tau))|_{n \in \set{N}}, \hat{\vect{u}}(\tau), \hat{\vect{a}}(\tau)$ be decision variables from any other policy and $\hat{r}_n(t) \triangleq r_n^{(\hat{f}_n(t))}(t),\hat{d}_n(t) \triangleq d_n^{(\hat{f}_n(t))}(t)$.  From \eqref{eq:square_drift} and definition \ref{def:quadratic_policy}, the drift-plus-penalty under quadratic policy is bounded by \eqref{eq:drift_under_quadratic} and is further bounded by \eqref{eq:drift_relaxed_quadratic} under any other policy as
\begin{align}
  & L(\queue(\tau+1)) - L(\queue(\tau)) - Vy_0(t)(\tau) \notag \\
  & \scriptstyle \leq \frac{1}{2} \sum_{n \in \set{N}} \bigl\{ \left[ K_n(\tau) - s_n^{(q)}(\tau) \right]^2 + \left[ K_n(\tau) - s_n^{(j)}(\tau) \right]^2 + \left[ K_n(\tau) + d_n(\tau) \right]^2  \notag \\
  & \scriptstyle \quad + \left[Q_n(\tau) - u_n(\tau) \right]^2 + \left[ Q_n(\tau) + s_n^{(q)}(\tau) \right]^2 + \sum_{m \in \set{N}} \left[Q_n(\tau) + a_{mn}(\tau) \right]^2  \notag \\
  & \scriptstyle \quad + \sum_{m \in \set{N}} \left[J_n(\tau) - a_{nm}(\tau) \right]^2 + [ J_n(\tau) + s_n^{(j)}(\tau) ]^2 - 2V r_n(\tau) + D_n(\tau)\bigr\} \label{eq:drift_under_quadratic} \\
  & \scriptstyle \leq \frac{1}{2} \sum_{n \in \set{N}} \bigl\{ \left[ K_n(\tau) - \hat{s}_n^{(q)}(\tau) \right]^2 + \left[ K_n(\tau) - \hat{s}_n^{(j)}(\tau) \right]^2 + \left[ K_n(\tau) + \hat{d}_n(\tau) \right]^2  \notag \\
  & \scriptstyle \quad + \left[Q_n(\tau) - \hat{u}_n(\tau) \right]^2 + \left[ Q_n(\tau) + \hat{s}_n^{(q)}(\tau) \right]^2 + \left[Q_n(\tau) + \sum_{m \in \set{N}} \hat{a}_{mn}(\tau) \right]^2  \notag \\
  & \scriptstyle \quad + \left[J_n(\tau) - \sum_{m \in \set{N}} \hat{a}_{nm}(\tau) \right]^2 + [ J_n(\tau) + \hat{s}_n^{(j)}(\tau) ]^2 - 2V \hat{r}_n(\tau) + D_n(\tau) \bigr\}. \label{eq:drift_relaxed_quadratic}
\end{align}

From the bounds \eqref{eq:linear_bound}, it follows that
\begin{align}
  & L(\queue(\tau+1)) - L(\queue(\tau)) - V y_0(\tau) \notag \\
  & \textstyle \leq \sum_{n \in \set{N}} \bigl\{ K_n(\tau) \left[ \hat{d}_n(\tau) - \hat{s}_n^{(q)}(\tau) - \hat{s}_n^{(j)}(\tau) \right] \notag \\
  & \text \quad + Q_n(\tau) \left[ \hat{s}_n^{(q)}(\tau) + \sum_{m \in \set{N}} \hat{a}_{mn}(\tau) - \hat{u}_n(\tau) \right] \notag \\
  & \textstyle \quad + J_n(\tau) \left[ \hat{s}_n^{(j)}(\tau) -\sum_{m \in \set{N}} \hat{a}_{nm}(\tau) \right] \notag \\
  & \textstyle \quad - V \hat{r}_n(\tau) \bigr\} + E \label{eq:drift_linear}
\end{align}
where
\begin{align}
  E
  & \triangleq \frac{1}{2} \mbox{$\sum_{n \in \set{N}}$} \Bigl\{ \prts{ s_n^{(q)\maxvar} + s_n^{(j)\maxvar} + d_n^{\maxvar} }^2 \notag \\ 
  & \quad + \prts{ s_n^{(q)\maxvar} + u_n^\maxvar + \mbox{$\sum_{m \in \set{N}}$} a_{mn}^\maxvar }^2 \notag \\
  & \quad + \prts{ s_n^{(j)\maxvar} + \mbox{$\sum_{m \in \set{N}}$} a_{nm}^\maxvar }^2 \Bigr\} \label{eq:E_const}
\end{align}

The derivations \eqref{eq:drift_under_quadratic}--\eqref{eq:drift_linear} show that applying the quadratic policy to the drift-plus-penalty expression leads to the bound \eqref{eq:drift_linear} which is valid for every other control policy.  However, the linear minimization of \eqref{eq:drift_linear}, which leads to the max-weight policy, does not resemble quadratic minimization of the actual drift-plus-penalty term \eqref{eq:pure_drift}.  The effects of the two policies are revealed in Sec. \ref{sec:simulation} where the quadratic policy leads to smaller queue backlogs.


As discussed in Sec. \ref{sec:model}, $\random(t)$ is i.i.d. over slots and is assumed further to have distribution $\pi(\random)$.  Define an $\random$-only policy as one that make a (possibly randomized) choice of decision variables based only on the observed $\random(t)$.  Then we customize an important theorem from \cite{Neely:SNObook}.

\begin{theorem}
  \label{thm:omega_only}
When problem \eqref{eq:optimization_problem} with stationary distribution $\pi(\random)$ is feasible, then for any $\delta > 0$ there exists an $\random$-only policy that chooses all controlled variables $(f_n^\ast(t), s_n^{(q)\ast}(t), s_n^{(j)\ast}(t))|_{n \in \set{N}}, \vect{u}^\ast(t), \vect{a}^\ast(t)$, and for all $n \in \set{N}$:
  \begin{align}
    \textstyle \expect{ y_0^\ast(t) } & \leq y_0^\optvar + \delta \label{eq:omega-1}\\
    \textstyle \expect{ d_n^\ast(t) - s_n^{(q)\ast}(t) - s_n^{(j)\ast}(t) } & \leq \delta \label{eq:omega-2}\\
    \textstyle \expect{ s_n^{(q)\ast}(t) + \sum_{m \in \set{N}} a_{mn}^\ast(t) - u_n^\ast(t) } & \leq \delta \label{eq:omega-3}\\
    \textstyle \expect{ s_n^{(j)\ast}(t) - \sum_{m \in \set{N}} a_{nm}^\ast(t) } & \leq \delta \label{eq:omega-4}
  \end{align}
where $y_0^\optvar$ is the optimal solution of problem \eqref{eq:optimization_problem}.  Also, $y_0^\ast(t) \triangleq \sum_{n \in \set{N}} r_n^\ast(t)$ when $r_n^\ast(t) \triangleq r_n^{(f_n^\ast(t))}(t)$ and $d_n^\ast(t) \triangleq d_n^{(f_n^\ast(t))}(t)$.
\end{theorem}

We additionally assume all constraints of the network can be achieved with $\epsilon$ slackness \cite{Neely:SNObook}:
\begin{assumption}
  \label{ass:slater}
  There are values $\epsilon > 0$ and $0 \leq y_0^{(\epsilon)} \leq y_0^\maxvar$ and an $\random$-only policy choosing all controlled variables $(f_n^\ast(t), s_n^{(q)\ast}(t), s_n^{(j)\ast}(t))|_{n \in \set{N}}, \vect{u}^\ast(t), \vect{a}^\ast(t)$ that satisfies for all $n \in \set{N}$:
  \begin{align}
    \textstyle \expect{ y_0^\ast(t) } & = y_0^{(\epsilon)} \label{eq:slater-1} \\
    \textstyle \expect{ d_n^\ast(t) - s_n^{(q)\ast}(t) - s_n^{(j)\ast}(t) } & \leq -\epsilon \label{eq:slater-2}\\
    \textstyle \expect{ s_n^{(q)\ast}(t) + \sum_{m \in \set{N}} a_{mn}^\ast(t) - u_n^\ast(t) } & \leq -\epsilon \label{eq:slater-3}\\
    \textstyle \expect{ s_n^{(j)\ast}(t) - \sum_{m \in \set{N}} a_{nm}^\ast(t) } & \leq -\epsilon \label{eq:slater-4}.
  \end{align}
\end{assumption}

\subsection{Performance Analysis}
Since our quadratic algorithm satisfies the bound \eqref{eq:drift_linear}, where the right-hand-side is in terms of any alternative policy $\left(\hat{h}_n(t), \hat{s}_n^{(q)}(t), \hat{s}_n^{(j)}(t) \right)|_{n \in \set{N}}, \hat{\vect{u}}(t), \hat{\vect{a}}(t)$, it holds for any $\random$-only policy $\left(h_n^\ast(t), s_n^{(q)\ast}(t), s_n^{(j)\ast}(t) \right)|_{n \in \set{N}}, \vect{u}^\ast(t), \vect{a}^\ast(t)$. Substituting an $\random$-only policy into \eqref{eq:drift_linear} and taking conditional expectations gives:
\begin{align}
  & \textstyle \expect{ L(\queue(\tau+1)) - L(\queue(\tau)) - V y_0(\tau) | \queue(\tau) } \\
  & \textstyle \leq \sum_{n \in \set{N}} \Bigl\{ K_n(\tau) \expect{\left. d_n^\ast(\tau) -s_n^{(q)\ast}(\tau) - s_n^{(j)\ast}(\tau) \right| \queue(\tau) } \notag \\
  & \textstyle \quad + Q_n(\tau) \expect{\left. s_n^{(q)\ast}(\tau) + \sum_{m \in \set{N}} a_{mn}^\ast(\tau) - u_n^\ast(\tau) \right| \queue(\tau) } \notag \\
  & \textstyle \quad + J_n(\tau) \expect{\left. s_n^{(j)\ast}(\tau) - \sum_{m \in \set{N}} a_{nm}^\ast(\tau) \right| \queue(\tau) } \notag \\
  & \textstyle \quad - V \expect{\left. r_n^\ast(\tau) \right| \queue(\tau) } \Bigr\} + E \notag \\
  & \textstyle \leq \sum_{n \in \set{N}} \Bigl\{ K_n(\tau) \expect{ d_n^\ast(\tau) -s_n^{(q)\ast}(\tau) - s_n^{(j)\ast}(\tau) } \notag \\
  & \textstyle \quad + Q_n(\tau) \expect{ s_n^{(q)\ast}(\tau) + \sum_{m \in \set{N}} a_{mn}^\ast(\tau) - u_n^\ast(\tau) } \notag \\
  & \textstyle \quad + J_n(\tau) \expect{ s_n^{(j)\ast}(\tau) - \sum_{m \in \set{N}} a_{nm}^\ast(\tau) } \notag\\
  & \textstyle \quad - V \expect{ r_n^\ast(\tau) } \Bigr\} + E \label{eq:dpp_w-only}
\end{align}
where we have used the fact that conditional expectations given $\queue(t)$ on the right-hand-side above are the same as unconditional expectations because $\random$-only policies do not depend on $\queue(t)$.

\begin{theorem}
  \label{thm:quad_bound}
  If Assumption \ref{ass:slater} holds, then the time-averaged total quality of information $\bar{y}_0$ is within $O(1/V)$ of optimality under the quadratic policy, while the total queue backlog grows with $O(V)$.
\end{theorem}

Theorem \ref{thm:quad_bound} is proven by substituting the $\random$-only policies from Theorem \ref{thm:omega_only} and Assumption \ref{ass:slater} into the right-hand-side of \eqref{eq:dpp_w-only}, as shown in the next subsections.

\subsubsection{Quality of Information vs. $V$}
\label{ssec:QoI_V}
Using the $\random$-only policy from \eqref{eq:omega-1}--\eqref{eq:omega-4} in the right-hand-side of \eqref{eq:dpp_w-only} gives:
\begin{align*}
  & \textstyle \expect{ L(\queue(\tau+1)) - L(\queue(\tau)) - V y_0(\tau) | \queue(\tau) } \\
  & \textstyle \leq E - V \left(y_0^{(\text{opt})} + \delta \right) + \delta \sum_{n \in \set{N}} \left[ K_n(\tau) + Q_n(\tau) + J_n(\tau) \right]
\end{align*}
This inequality is valid for every $\delta > 0$.  Therefore
\begin{equation*}
  \textstyle \expect{ L(\queue(\tau+1)) - L(\queue(\tau)) - V y_0(\tau) | \queue(\tau) } \leq E - V y_0^\optvar.
\end{equation*}
Taking an expectation and summing from $\tau = 0$ to $t-1$:
\begin{equation*}
  \textstyle \expect{ L(\queue(t)) - L(\queue(0)) - V \sum_{\tau = 0}^{t-1} y_0(\tau) } \leq E t - Vt y_0^\optvar.
\end{equation*}
With rearrangement and $L(\queue(t)) \geq 0$, it follows that
\begin{equation*}
  \sum_{\tau = 0}^{t-1} \expect{y_0(\tau)} \geq -\frac{E t}{V} + t y_0^\optvar - \frac{L(\queue(0))}{V}.
\end{equation*}

Dividing by $t$ and taking limit as $t$ approaches infinity, the performance of the quadratic policy is lower bounded by
\begin{equation}
  \label{eq:performance_qoi}
  \liminf_{t \rightarrow \infty} \frac{1}{t} \sum_{\tau = 0}^{t-1} \expect{ y_0(\tau) } \geq -\frac{ E }{ V } + y_0^\optvar.
\end{equation}

This shows that the system can be pushed to the optimality $y_0^\optvar$ by increasing $V$ under the quadratic policy.

\subsubsection{Total Queue Backlog vs. $V$}
\label{ssec:backlog_V}
Now consider the existence of an $\random$-only policy with Assumption \ref{ass:slater}.  Using \eqref{eq:slater-1}--\eqref{eq:slater-4} in the right-hand-side of \eqref{eq:dpp_w-only} gives:
\begin{multline*}
  \expect{ L(\queue(\tau+1)) - L(\queue(\tau)) - V y_0(\tau) | \queue(\tau) } \\
  \textstyle \leq E - V y_0^{(\epsilon)} - \epsilon \sum_{n \in \set{N}} \left[ K_n(\tau) + Q_n(\tau) + J_n(\tau) \right].
\end{multline*}

Taking expectation and summing from $\tau = 0$ to $t-1$:
\begin{multline*}
  \expect{ L(\queue(t)) - L(\queue(0)) - V  \mbox{$\sum_{\tau = 0}^{t-1}$} y_0(\tau)} \\
  \textstyle \leq E t - Vt y_0^{(\epsilon)} - \epsilon \sum_{\tau = 0}^{t-1} \sum_{n \in \set{N}} \expect{ K_n(\tau) + Q_n(\tau) + J_n(\tau) }
\end{multline*}

With rearrangement and $L(\queue(t)) \geq 0$, it follows that
\begin{align*}
  & \textstyle \sum_{\tau = 0}^{t-1} \sum_{n \in \set{N}} \expect{K_n(\tau) + Q_n(\tau) + J_n(\tau)} \\
  & \textstyle \quad \leq \frac{E t}{\epsilon} + \frac{V}{\epsilon}\left( \sum_{\tau = 0}^{t-1} \expect{y_0(\tau)} - t y_0^{(\epsilon)} \right) + \frac{\expect{L(\queue(0))}}{\epsilon} \\
  & \textstyle \quad \leq \frac{E t}{\epsilon} + \frac{V}{\epsilon}\left( t y_0^\maxvar - t y_0^{(\epsilon)} \right) + \frac{\expect{L(\queue(0))}}{\epsilon}.
\end{align*}

Dividing by $t$ and taking limit as $t$ approaches infinity, the time-averaged total queue backlog is bounded by
\begin{multline}
  \label{eq:performance_queue}
  \textstyle \lim \sup_{t \rightarrow \infty} \frac{1}{t} \sum_{\tau = 0}^{t-1} \sum_{n \in \set{N}} \expect{ K_n(\tau) + Q_n(\tau) + J_n(\tau) } \\
  \leq \frac{E}{\epsilon} + \frac{V}{\epsilon} \left( y_0^\maxvar - y_0^{(\epsilon)} \right).
\end{multline}

This shows that the overall queue length tends to increase linearly as $V$ is increased.  This is an asymptotic bound which shows that every queue is strongly stable, and the network is strongly stable.

The $V$ parameter in \eqref{eq:performance_qoi} and \eqref{eq:performance_queue} affects the performance trade-off $[O(1/V), O(V)]$ between quality of information and total queue backlog.  These results are similar to those that can be derived under the max-weight algorithm.  However, simulation in the next section shows significant reduction of queue backlog under the quadratic policy.

\subsection{Deterministic bounds of queue lengths}
Here we show that, in addition to the average queue size bounds derived in the previous subsection, our algorithm also yields deterministic worst-case queue size bounds which is summarized in the following lemma.  Define $K_n^\maxvar = \max_{f \in \set{F}} \frac{2 V r_n^{(f)} - d_n^{(f)2}}{2 d_n^{(f)}} + d_n^\maxvar$ for $n \in \set{N}$, and $Q_n^\maxvar \defequiv \max\prts{ K_n^\maxvar, \prtc{K_m^\maxvar}_{m \in \set{N}}} + \sum_{m \in \set{N}} a_{mn}^\maxvar + s_n^\maxvar$.

\begin{lemma}
  For all devices $n \in \set{N}$ and all slots $t \geq 0$, we have:
  \begin{align}
    K_n(t) & \leq K_n^\maxvar \label{eq:Kbound}\\
    J_n(t) & \leq K_n^\maxvar \label{eq:Jbound}\\
    Q_n(t) & \leq Q_n^\maxvar \label{eq:Qbound}
  \end{align}
provided that these inequalities hold at $t = 0$.
\end{lemma}
\begin{IEEEproof}
  The bounds \eqref{eq:Kbound}--\eqref{eq:Qbound} are proved in Section \ref{sssec:Kbound}--\ref{sssec:Qbound} respectively.
\end{IEEEproof}

\subsubsection{Input Queue}
\label{sssec:Kbound}
From the admission-control problem \eqref{eq:admission_control}, if $\prtr{ r_n(t), d_n(t) } = (0, 0)$, then the objective value of the problem is $K_n(t)^2$.  Therefore, device $n$ only chooses $\prtr{ r_n(t), d_n(t) } \neq (0, 0)$ when 
\begin{align*}
  \prts{ K_n(t) + d_n(t) }^2 - 2 V r_n(t) & \leq K_n(t)^2 \notag \\
  2 K_n(t) d_n(t) + d_n(t)^2 - 2 V r_n(t) & \leq 0 \notag
\end{align*}
{\vskip -1.7em}
\begin{align}
  K_n(t) & \leq \frac{2 V r_n(t) - d_n(t)^2}{2 d_n(t)} \notag \\
         & \leq \max_{f \in \set{F}} \frac{2 V r_n^{(f)} - d_n^{(f)2}}{2 d_n^{(f)}} \label{eq:Kmax}
\end{align}
This implies that device $n$ can only obtain data when \eqref{eq:Kmax} holds, and receives no new data otherwise.

Fix $t$, and assume $K_n(t) \leq K_n^\maxvar$ for this slot $t$.  From \eqref{eq:input_queue}, there are two cases to consider.\\
i) If $0 \leq K_n(t) \leq K_n^\maxvar - d_n^\maxvar$ then \eqref{eq:Kmax} holds and $K_n(t+1) = K_n(t) + d_n(t) \leq K_n^\maxvar$.\\
ii) If $K_n(t) > K_n^\maxvar - d_n^\maxvar$, then \eqref{eq:Kmax} does not hold and $K_n(t+1) = K_n(t) \leq K_n^\maxvar$. Thus, given that $K_n(0) \leq K_n^\maxvar$, $K_n(t) \leq K_n^\maxvar$ for all $t \geq 0$ by mathematical induction.

\subsubsection{Relay Queue}
\label{sssec:Jbound}
Fix $t$ and assume for each device $n \in \set{N}$ that $J_n(t) \leq K_n^\maxvar$ for this slot $t$.  From the closed form solution \eqref{eq:closed_form_snj} and \eqref{eq:relay_queue}, there are three cases to consider.\\
i) When $K_n(t) - J_n(t) \leq 0$, then $s_n^{(j)}(t) = 0$, and 
\begin{align*}
  J_n(t+1) & \leq \max \prts{ J_n(t) + s_n^{(j)}(t), 0} \\
  & = J_n(t) \leq K_n^\maxvar.
\end{align*}
ii) When $K_n(t) - J_n(t) \geq 2 s_n^{(j)\maxvar}$ (or $J_n(t) \leq K_n(t) - 2 s_n^{(j)\maxvar}$), then $s_n^{(j)}(t) = s_n^{(j)\maxvar}$, and
\begin{align*}
  J_n(t+1) & \leq \max\prts{ J_n(t) + s_n^{(j)}(t), 0 } \\
  & \leq \max\prts{ K_n(t) - s_n^{(j)\maxvar}, 0} \\
  & \leq K_n(t) \leq K_n^\maxvar.
\end{align*}
iii) When $0 < K_n(t) - J_n(t) < 2 s_n^{(j)\maxvar}$, then $s_n^{(j)}(t) \leq \left\lceil \frac{K_n(t) - J_n(t)}{2} \right\rceil$, and
\begin{align*}
  J_n(t+1) & \leq \max\prts{ J_n(t) + s_n^{(j)}(t), 0 } \\
  & \leq \max\prts{ \left\lceil \frac{K_n(t) + J_n(t)}{2} \right\rceil , 0 } \\
  & \leq K_n(t) \leq K_n^\maxvar.
\end{align*}
Thus, given that $J_n(0) \leq K_n^\maxvar$, $J_n(t) \leq K_n^\maxvar$ for all $t \geq 0$ by mathematical induction.

\subsubsection{Uplink Queue}
\label{sssec:Qbound}
To provide a general upper bound for the uplink queue, we assume that all relay channels are orthogonal.  This implies every device $n \in \set{N}$ can transmit and receive relayed data simultanously. 

Fix $t$ and assume $Q_n(t) \leq Q_n^\maxvar$ for this slot $t$.  Then consider $Q_n(t+1)$ from \eqref{eq:uplink_queue}.\\
i) When $Q_n(t) \geq \max\prts{ K_n^\maxvar, \prtc{K_m^\maxvar}_{m \in \set{N}}}$, from \eqref{eq:closed_form_snq} and \eqref{eq:closed_form_anmt}, it follows that $s_n^{(q)}(t) = 0$ and $a_{mn}(t) = 0$ for all $m \in \set{N}$, so $Q_n(t+1) \leq Q(t) \leq Q_n^\maxvar$.\\
ii) When $Q_n(t) < \max\prts{ K_n^\maxvar, \prtc{K_m^\maxvar}_{m \in \set{N}}}$, then this queue may received data $s_n^{(q)}(t)$ and $a_{mn}(t)$ for some $m \in \set{N}$, so
\begin{align*}
  Q_n(t+1) & \leq \max\prts{ Q_n(t) + s_n^{(q)}(t), 0 } + \sum_{m \in \set{N}} a_{mn}(t)\\
  & \leq Q(t) + s_n^{(q)\maxvar} + \sum_{m \in \set{N}} a_{mn}^\maxvar\\
  & \leq Q_n^\maxvar.
\end{align*}
Thus, given $Q_n(0) \leq Q_n^\maxvar$,  $Q_n(t) \leq Q_n^\maxvar$ for all $t \geq 0$ by mathematical induction.

\section{Simulation}
\label{sec:simulation}
\begin{figure}
  \centering
  \includegraphics[scale=0.9]{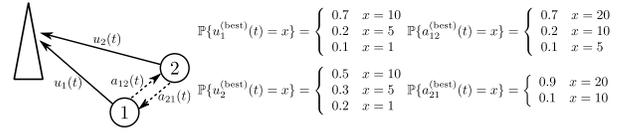}
  \caption{Small network with independent channels with distributions shown.}
  \label{fig:sim_network}
\end{figure}

Simulation under the proposed quadratic policy and the standard max-weight policy is performed over a small network in Fig. \ref{fig:sim_network}.  The network contains two devicess, $\set{N} = \{ 1, 2 \}$.  Each device has the other as its neighbor, so $\set{H}_1 = \{ 2 \}$ and $\set{H}_2 = \{ 1 \}$.  An event occurs in every slot with probability $\theta = 0.3$.  We assume all uplink and relay channels are orthogonal.  The uplink channel distribution for device $1$ is better than that of device $2$ as in Fig. \ref{fig:sim_network}.

The constraints are $u_n(t) \in \{0, \ldots, u_n^\bestvar(\vect{\eta}(t))\}$ for $n \in \set{N}$.  Also, $a_{12}(t) \in \{0, \ldots, a_{12}^\bestvar(\vect{\eta}(t))\}$ and $a_{21}(t) \in \{0, \ldots, a_{21}^\bestvar(\vect{\eta}(t))\}$.  Then set $s_n^{(q)\maxvar} = s_n^{(j)\maxvar} = 30$.  The feasible set of formats is $\set{F} = \{0, 1, 2, 3\}$ with constant options given by $(d_n^{(0)}, r_n^{(0)}) = (0, 0), (d_n^{(1)}, r_n^{(1)}) = (100, 20), (d_n^{(2)}, r_n^{(2)}) = (50, 15), (d_n^{(3)}, r_n^{(3)}) = (10, 10)$ whenever there is an event. 

\begin{figure}
  \centering
  \includegraphics[scale = 0.42]{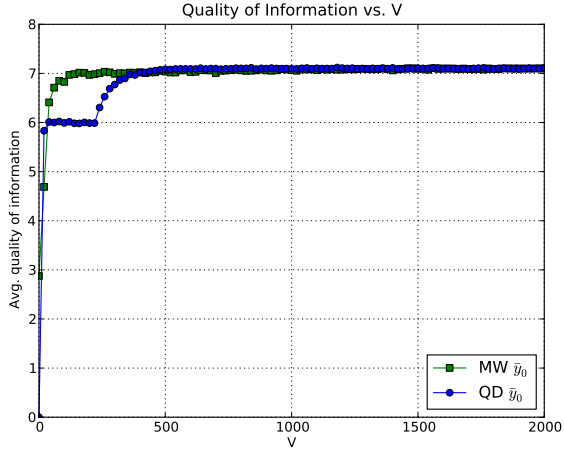}
  \caption{Quality of Information versus $V$ under the quadratic (QD) and max-weight (MW) policies}
  \label{fig:QoI_plot}
\end{figure}

The simulation is performed according to the algorithm in Sec. \ref{ssec:algorithm}.  The time-averaged quality of information under the quadratic and max-weight policies are shown in Fig. \ref{fig:QoI_plot}.  From the plot, the values of $\bar{y}_0$ under both policies converge to optimality following the $O(1/V)$ performance bound.

\begin{figure}
  \centering
  \includegraphics[scale = 0.42]{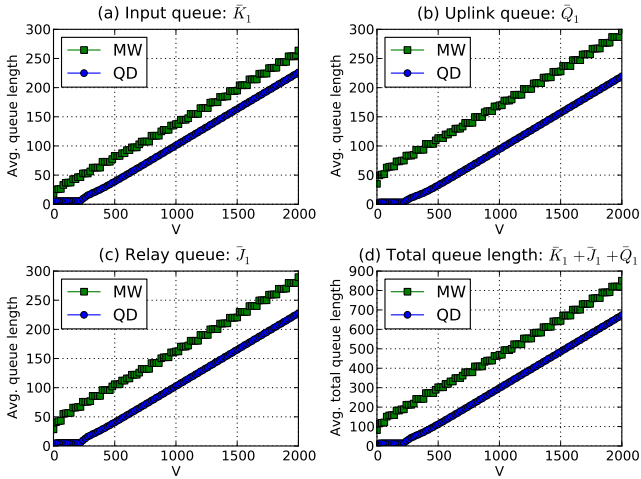}
  \caption{Averaged backlog in queues versus $V$ and system quality versus backlog under the quadratic (QD) and max-weight (MW) policies}
  \label{fig:all_backlog_plot}
\end{figure}

Fig. \ref{fig:all_backlog_plot}abc reveals queue lengths in the input, uplink, and relay queues of device 1 under the quadratic and max-weight policies.  At the same $V$, the quadratic policy reduces queue lengths by a significant constant compared to the cases under the max-weight policy.  The plot also shows the growth of queue lengths with parameter $V$, which follows the $O(V)$ bound of the queue length.  Fig. \ref{fig:all_backlog_plot}d shows the average total queue length in device 1 under the quadratic and max-weight policies.

\begin{figure}
  \centering
  \includegraphics[scale = 0.42]{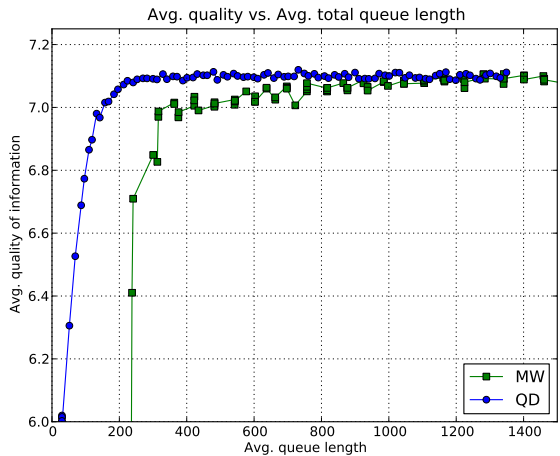}
  \caption{The system obtains average quality of information while having average total queue length}
  \label{fig:all_QoI_V}
\end{figure}

Fig. \ref{fig:all_QoI_V} shows that the quadratic policy can achieve near optimality with significantly smaller total system backlog compared to the case under the max-weight policy.  This shows a significant advantage, which in turn affects memory size and packet delay. 

\begin{figure}
  \centering
  \includegraphics[scale = 0.9]{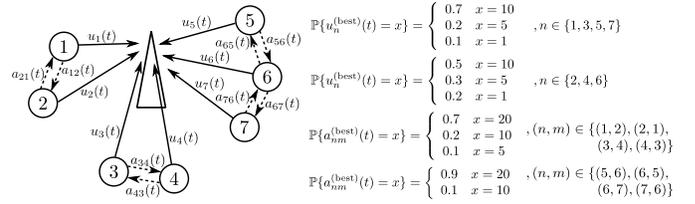}
  \caption{Larger network with independent channels with distributions shown}
  \label{fig:medium_sce}
\end{figure}

\begin{figure}
  \centering
  \includegraphics[scale = 0.42]{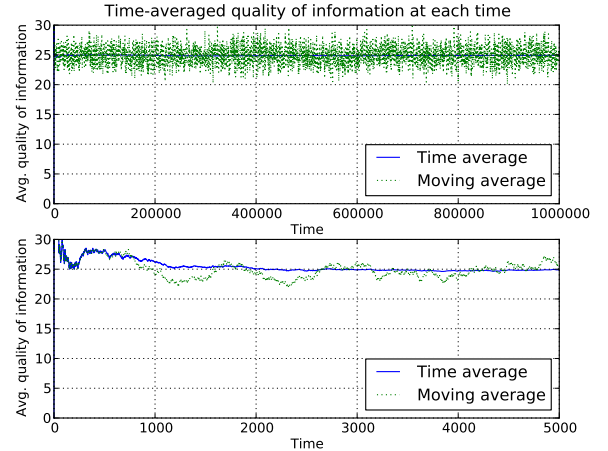}
  \caption{Convergence of time-averaged quality of information.  The interval of the moving average is 500 slots.}
  \label{fig:medium_rewards}
\end{figure}

Another larger network shown in Fig. \ref{fig:medium_sce} is simulated to observe convergence of the proposed algorithm.  As in the small network scenario, the same probability of event occurrence $\theta = 0.3$ is set.  Channel distributions are configured in Fig. \ref{fig:medium_sce}.  For $V = 800$, the time-averaged quality of information is $25.00$ after $10^6$ time slots as shown in the upper plot of Fig. \ref{fig:medium_rewards}.  The lower plot in Fig. \ref{fig:medium_rewards} illustrates the early period of the simulation to illustrate convergence time.

\section{Linear Programs by Quadratic Policy}
\label{sec:linear_program}
The generality of the quadratic policy is illustrated in this section. The policy is applied to solve linear programs which is one application of the Lyapunov optimization \cite{Neely:NOW}.

\subsection{Problem Transformation}
The following static linear program is considered where $\prtr{x_i}_{i=1}^n$ are decision variables and $\prtr{a_{ji}}|_{j=1,i=1}^{j=m,i=n}$, $\prtr{b_j}|_{j=1}^m$, $\prtr{c_i}|_{i=1}^n$, $\prtr{x_i^\maxvar}|_{i=1}^n$ are constants.
\begin{align}
  \maximize \quad 
  & \sum_{i = 1}^n c_i x_i \label{eq:linearform}\\
  \subjectto \quad
  & \sum_{i = 1}^n a_{ji} x_i \leq b_j, \quad\quad j \in \prtc{1, \dotsc, m} \notag\\
  & 0 \leq x_i \leq x_i^\maxvar, \quad\quad i \in \{1, \dotsc, n\} \notag
\end{align}

In order to solve \eqref{eq:linearform}, the following time-averaged optimization problem is solved by using the Lyapunov optimization technique.
\begin{align}
  \maximize \quad
  & \bar{y}_0 = \sum_{i = 1}^n c_i \bar{x}_i \label{eq:linear_avg_problem}\\
  \subjectto \quad
  & \sum_{i = 1}^n a_{ji} \bar{x}_i \leq b_j, \quad\quad j \in \prtc{1, \dotsc, m} \notag\\
  & 0 \leq x_i(t) \leq x_i^\maxvar, \quad\quad i \in \prtc{1, \dotsc, n}, t \geq 0 \notag
\end{align}

Solutions from the static problem \eqref{eq:linearform} and the time-averaged problem \eqref{eq:linear_avg_problem} are equivalent because using a solution $x_i$ to the static problem for every $t$ in the time-averaged problem leads to $\bar{x}_i = x_i$ for $i \in \prtc{1, \dotsc, n}$ and every constraint in the time-averaged problem is satisfied.  The time-averaged objective function is also maximized, since the time average of the linear function is equal to the function of the time averages.  On the other hand, a solution to the time-averaged problem is a solution of the static problem because it satisfies all constraints and maximizes the same objective function.

To solve problem \eqref{eq:linear_avg_problem}, a concept of virtual queue is used \cite{Neely:NOW}.  Let $x_i(t)$ be chosen every slot $t$ in the interval $0 \leq x_i(t) \leq x_i^\maxvar$, and define $\bar{x}_i(t)$ for $i \in \prtc{1, \dotsc, n}$ as
\begin{equation*}
  \bar{x}_i(t) \defequiv \frac{1}{t} \sum_{\tau = 0}^{t-1} x_i(\tau), \quad i \in \prtc{1, \dotsc, n}.
\end{equation*}

Define virtual queue
\begin{align}
  Z_j(t+1)
  & = \max\prts{ Z_j(t) + \sum_{i = 1}^n a_{ji} x_i(t) - b_j , 0}  \label{eq:virtual_con1}
\end{align}
for all $j \in \prtc{1, \dotsc, m}$.  It follows that
\begin{align*}
  Z_j(\tau+1) & = \max \prts{ Z_j(\tau) + \sum_{i = 1}^n a_{ji} x_i(\tau) - b_j , 0} \\
              & \geq Z_j(\tau) + \sum_{i=1}^n a_{ji} x_i(\tau) - b_j \\
  Z_j(\tau+1) - Z_j(\tau) & \geq \sum_{i=1}^n a_{ji} x_i(\tau) - b_j.
\end{align*}
Summing from $\tau = 0$ to $t-1$, and dividing by $t$:
\begin{align}
  Z_j(t) - Z_j(0)    & \geq \sum_{\tau=0}^{t-1} \sum_{i=1}^n a_{ji} x_i(\tau) - t b_j \notag\\
  \frac{Z_j(t)}{t}   & \geq \frac{1}{t} \sum_{\tau=0}^{t-1} \sum_{i=1}^n a_{ji} x_i(\tau) - b_j \notag \\
                     & = \sum_{i=1}^n a_{ji} \bar{x}_i(t) - b_j \label{eq:rate_stable_bound},
\end{align}
where we assume that $Z_j(0) \geq 0$. It follows that if $\frac{Z(t)}{t} \rightarrow 0$ (so that each queue is ``rate stable''), the desired time-average inequality constraint is satisfied. 

Then let $\queue(t) = \prtr{ Z_j(t) }|_{j=1}^m$ be a vector of all virtual queues and
\begin{equation*}
  y_0(t) = \sum_{i=1}^n c_i x_i(t)
\end{equation*}
be the objective function whose time average is to be minimized according to the problem \eqref{eq:linear_avg_problem}.  Define a time-averaged objective value up to iteration $t$ by
\begin{align*}
  \bar{y}_0(t) & \defequiv \frac{1}{t} \sum_{\tau = 0}^{t-1} y_0(\tau).
\end{align*}
Similar to Sec. \ref{ssec:stochastic_optimization}, let $\bar{y}_0 \triangleq \lim_{t \rightarrow \infty} \bar{y}_0(t)$ be and $\bar{x}_i \triangleq \lim_{t \rightarrow \infty} \bar{x}_i(t)$ be their asymptotic averages.

\subsection{Lyapunov Optimization}
To solve \eqref{eq:linear_avg_problem}, the drift-plus-penalty for this problem is bounded by Lemma \ref{lem:absbounds} as
\begin{align}
  & L(\queue(t+1)) - L(\queue(t)) - V y_0(t) \notag\\
  & \textstyle \quad = \frac{1}{2} \sum_{j = 1}^m \prtc{ Z_j^2(t+1) - Z_j^2(t) - 2 V c_i x_i(t) } \notag\\
  & \textstyle \quad = \frac{1}{2} \sum_{j = 1}^m \bigl\{ \max \prts{ Z_j(t) + \sum_{i=1}^n a_{ji} x_i(t) - b_j , 0 }  \notag\\
  & \textstyle \quad\quad - Z_j^2(t) - 2 V c_i x_i(t) \bigr\} \label{eq:linearpure}\\
  & \textstyle \quad \leq \frac{1}{2} \sum_{j=1}^m \bigl\{ \sum_{i=1}^n \prts{Z_i(t) + a_{ji} x_i(t)}^2 + \prts{Z_j(t) - b_j}^2 \notag\\
  & \textstyle \quad\quad - (n+1)Z_j^2(t) - 2 V c_i x_i + H_j \bigr\} \notag\\
  & \textstyle \quad = \frac{1}{2} \sum_{i=1}^n \prtc{ \sum_{j=1}^m \prts{ Z_j(t) + a_{ji} x_i(t) }^2 - 2 V c_i x_i(t) } \notag\\
  & \textstyle \quad\quad + \frac{1}{2} \sum_{j=1}^m \prtc{ \prts{ Z_j(t) - b_j }^2 - (n+1) Z_j^2(t) + H_j }, \label{eq:linearqd}
\end{align}
where
\begin{multline*}
  H_j = 2 \bigl\{ \mbox{$\sum_{i=1}^n \sum_{i'=1}^{i-1}$} \abs{a_{ji}} \abs{a_{ji'}} x_i^\maxvar x_{i'}^\maxvar \\
  + \mbox{$\sum_{i=1}^n$} \abs{a_{ji}} \abs{b_j} x_i^\maxvar \bigr\}
\end{multline*}
for $j \in \prtc{1, \dotsc, m}$.  From \eqref{eq:linearqd}, the quadratic policy minimize the drift-plus-penalty every iteration, and this minimization is
\begin{align}
  \minimize \quad
  & \textstyle \sum_{i=1}^n \prtc{ \sum_{j=1}^m \prts{ Z_j(t) + a_{ji} x_i(t) }^2 - 2 V c_i x_i(t) } \label{eq:linear_qd_problem}\\
  \subjectto \quad  & 0 \leq x_i(t) \leq x_i^\maxvar \quad\quad i \in \prtc{1, \dotsc, n}. \notag
\end{align}


Again, because problem's structure and the fully separable property of the quadratic policy, problem \eqref{eq:linear_qd_problem} can be solved separately for each $x_i(t)$. A closed form solution of each $x_i(t)$ for $i \in \prtc{1, \dotsc, n}$ is
\begin{align*}
  x_i(t) = \max\prts{ \min\prts{ \frac{c_i V - \sum_{j=1}^m a_{ji} Z_j(t)}{\sum_{j=1}^m a_{ji}^2}, x_i^\maxvar}, 0}.
\end{align*}

\subsection{Algorithm}
An algorithm to solve problem \eqref{eq:linear_avg_problem}, which also solves \eqref{eq:linearform}, is the following.
\begin{algorithm}
  Initialize $\prtc{\queue(0)} = \vect{0}$ \\
  $t = 0$ \\
  \ForEach{iteration $t \geq 0$} {
    \tcp{Update decision variables}
    \ForEach{$i \in \prtc{1, \dotsc, n}$} {
      $\textstyle   x_i(t) = \max\prts{ \min\prts{ \frac{c_i V - \sum_{j=1}^m a_{ji} Z_j(t)}{\sum_{j=1}^m a_{ji}^2}, x_i^\maxvar}, 0}$ \\
      }
    \tcp{Update virtual queues}
    \ForEach{$j \in \prtc{1, \dotsc, m}$} {
      $\textstyle Z_j(t+1) = \max\prts{ Z_j(t) + \prtr{ \sum_{i=1}^n a_{ji} x_i(t) - b_j }, 0 }$ \\
      }
    $t = t+1$ \\
  }
  \caption{Linear programming by quadratic policy}
  \label{alg:linear}
\end{algorithm}

\subsection{Convergence Analysis}
Since our policy chooses $x_i(t) \in [0, x^\maxvar]$ every slot to minimize the right-hand-side of \eqref{eq:linearqd}, this right-hand-side is less than or equal to the corresponding value with any other feasible decision $x_i^\ast \in [0, x^\maxvar]$:
\begin{align}
  & L(\queue(\tau+1)) - L(\queue(\tau)) - V y_0(\tau) \notag\\
  & \textstyle \leq \frac{1}{2} \sum_{i=1}^n \prtc{ \sum_{j=1}^m \prts{ Z_j(t) + a_{ji} x_i^\ast(t) }^2 - 2 V c_i x_i^\ast(t) } \notag\\
  & \textstyle \quad + \frac{1}{2} \sum_{j=1}^m \prtc{ \prts{ Z_j(t) - b_j }^2 - (n+1) Z_j^2(t) + H_j } \notag\\
  & \textstyle \leq \sum_{j=1}^m Z_j(\tau) \prts{ \sum_{i = 1}^n a_{ji} x_i^\ast(\tau) - b_j } - V \sum_{i = 1}^n c_i x_i^\ast(\tau) + E \label{eq:linearmw} 
\end{align}
where
\begin{align*}
  E = \sum_{j=1}^m \prts{ \sum_{i=1}^n \abs{a_{ji}} x_i^\maxvar + b_j }^2
\end{align*}
and the final inequality uses \eqref{eq:linear_bound}.

Assume that problem \eqref{eq:linearform} has $\vect{x}^\ast = \prtr{ x_i^\ast }|_{i=1}^n$ as an optimal solution and $y_0^\optvar$ as the optimal cost.  This optimal solution has the following properties:
\begin{align*}
  y_0^\optvar & = \sum_{i=1}^n c_i x_i^\ast \\
  \sum_{i=1}^n a_{ji} x_i^\ast & \leq b_j, \quad\quad j \in \prtc{1, \dotsc, m}.
\end{align*}

By applying $x_i(t) = x_i^\ast$ every iteration, the bound \eqref{eq:linearmw} becomes
\begin{align*}
  & L(\queue(\tau+1)) - L(\queue(\tau)) - V y_0(\tau) \notag\\
  & \quad \leq -V y_0^\optvar + E \notag
\end{align*}

Summing from $\tau = 0$ to $t-1$ and rearranging lead to
\begin{align*}
  L(\queue(t)) - L(\queue(0)) - V\sum_{\tau=0}^{t-1} y_0(\tau) & \leq Et - Vt y_0^\optvar
\end{align*}
and
\begin{equation*}
  \sum_{\tau=0}^{t-1} y_0(\tau) \geq \frac{L(\queue(t)) - L(\queue(0)) - Et}{V} + t y_0^\optvar.
\end{equation*}
Since Algorithm \ref{alg:linear} initializes $\queue(0) = \vect{0}$, $L(\queue(0)) = 0$ and also $L(\queue(t)) \geq 0$. Then, dividing by $t$ leads to
\begin{align}
  \frac{1}{t} \sum_{\tau=0}^{t-1} y_0(\tau) & \geq \frac{L(\queue(t))}{tV} - \frac{E}{V} + y_0^\optvar \label{eq:linear_cost_bound_withL} \\
  \frac{1}{t} \sum_{\tau=0}^{t-1} y_0(\tau) & \geq - \frac{E}{V} + y_0^\optvar. \label{eq:linear_cost_bound}
\end{align}
The bound \eqref{eq:linear_cost_bound} shows that, when $V$ is large, the time-averaged objective value from the algorithm approaches the optimal objective value.

Since the feasible set of $\prtr{x_i(t)}_{i=1}^n$ in problem \eqref{eq:linear_avg_problem} is bounded, there exist some $y_0^\maxvar \geq 0$ such that $y_0(t) \leq y_0^\maxvar$ for all $t$.  Then, the bound \eqref{eq:linear_cost_bound_withL} can also be rearranged to be
\begin{align}
  L(\queue(t))
  & \leq Et + V \sum_{\tau=0}^{t-1} y_0(\tau) - Vt y_0^\optvar \notag\\
  \sum_{j = 1}^m Z_j^2(t)
  & \leq 2Et + 2Vt\prts{ y_0^\maxvar - y_0^\optvar } \notag\\
  Z_j(t)
  & \leq \sqrt{ 2Et + 2Vt\prts{ y_0^\maxvar - y_0^\optvar } } \notag\\
  \frac{Z_j(t)}{t}
  & \leq \sqrt{ \frac{1}{t} \prtc{ 2E + 2V\prts{ y_0^\maxvar - y_0^\optvar }} }. \notag
\end{align}
From \eqref{eq:rate_stable_bound}, it follows that
\begin{equation}
  \sum_{i = 1}^n a_{ji} \bar{x}_i(t) - b_j \leq \sqrt{ \frac{1}{t} \prtc{2E + 2V\prts{ y_0^\maxvar - y_0^\optvar }} }. \label{eq:linear_cons_bound}
\end{equation}
The bound \eqref{eq:linear_cons_bound} shows that the constraints of problem \eqref{eq:linearform} are asymptotically satisfied as $t$ approaches infinity.

When the number of iterations is limited, we can obtain convergence results in this case by assuming $V = 1/\varepsilon$ and $t = 1/\varepsilon^3$ and consider \eqref{eq:linear_cost_bound} and \eqref{eq:linear_cons_bound}.  This leads to
\begin{align*}
  y_0^\optvar - \frac{1}{t} \sum_{\tau=0}^{t-1} y_0(\tau) \leq E\varepsilon = O(\varepsilon). 
\end{align*}
and
\begin{align}
  \sum_{i = 1}^n a_{ji} \bar{x}_i(t) - b_j & \leq \sqrt{ \frac{1}{1/\varepsilon^3} \prtc{2E + 2/\varepsilon \times \prts{ y_0^\maxvar - y_0^\optvar }}} \notag \\
  & = O(\varepsilon). \label{eq:linear_cons_bound_limited}
\end{align}
Therefore, using $O(1/\varepsilon^3)$ iterations ensures the time-averaged value of $\bar{y}_0(t)$ is within $O(\varepsilon)$ of the optimal value $y_0^\optvar$, and all constraints are within $O(\varepsilon)$ of being satisfied.  However, This $O(1/\varepsilon^3)$ tradeoff can be improved to an $O(1/\varepsilon^2)$ tradeoff if the problem \eqref{eq:linear_avg_problem} satisfies a mild ``Slater assumption'' as the following.

\begin{assumption}
  \label{ass:staticslater}
  There are values $\epsilon > 0$ and $y_0^\optvar \leq y_0^{(\epsilon)} \leq y_0^\maxvar$ and a static policy choosing $( x_i^\ast )_{i=1}^n$ every iteration that satisfies:
  \begin{align*}
    \textstyle y_0^\ast(t)                        & = y_0^{(\epsilon)} \\
    \textstyle \sum_{i = 1}^n a_{ji} x_i^\ast - b_j & \leq -\epsilon \quad\quad j \in \prtc{1, \dotsc, m} \\
    \textstyle 0 \leq x_i^\ast                    & \leq x_i^\maxvar. 
  \end{align*}
\end{assumption}
In fact, this assumption is a static version of Assumption \ref{ass:slater} and is similar to a general Slater condition in the convex optimization theory \cite{Bertsekas:Convex}.

Applying Assumption \eqref{ass:staticslater} to \eqref{eq:linearmw}, it follows that
\begin{align*}
  & L(\queue(\tau+1)) - L(\queue(\tau)) - V y_0(\tau) \notag\\
  & \textstyle \quad \leq - V \sum_{i = 1}^n c_i x_i^\ast + \sum_{j=1}^m Z_j(\tau) \prts{ \sum_{i = 1}^n a_{ji} x_i^\ast - b_j } + E \\
  & \textstyle \quad \leq - V y_0^{(\epsilon)} - \epsilon \sum_{j = 1}^m Z_j(\tau) + E
\end{align*}

From triangle inequality, $\lVert \vect{Z}(\tau) \rVert \leq \sum_{j=1}^m Z_j$, the above inequality is
\begin{align*}
  & L(\queue(\tau+1)) - L(\queue(\tau)) - V y_0(\tau) \notag\\
  & \textstyle \quad \leq - V y_0^{(\epsilon)} - \epsilon \lVert \vect{Z}(\tau) \rVert + E
\end{align*}

Since $L(\queue(\tau)) = \frac{1}{2}\lVert \vect{Z}(\tau) \rVert^2$, we have:
\begin{align*}
  \lVert \vect{Z}(\tau+1) \rVert^2 - \lVert \vect{Z}(\tau) \rVert^2 \leq 2\prts{V \prtr{y_0(\tau) - y_0^{(\epsilon)}} - \epsilon \lVert \vect{Z}(\tau) \rVert + E}.
\end{align*}

If $\lVert \vect{Z}(\tau) \rVert \geq \frac{ V \prtr{ y_0^\maxvar - y_0^{(\epsilon)} } + E}{\epsilon}$, then 
\begin{align*}
  \lVert \vect{Z}(\tau+1) \rVert^2 - \lVert \vect{Z}(\tau) \rVert^2 & \leq 0.
\end{align*}
Since $\sum_{j=1}^m Z_j^2(\tau) = \lVert \vect{Z}(\tau) \rVert^2$, the above inequality implies that the value of $\sum_{j = 1}^m Z_j^2(\tau)$ is not increased in the next iteration.  Therefore, the value of each $Z_j(\tau)$ is bounded by, for all $\tau \geq 0$,
\begin{equation*}
  Z_j(\tau) \leq \frac{ V \prtr{ y_0^\maxvar - y_0^{(\epsilon)}} + E}{\epsilon} + \sum_{i = 1}^n |a_{ji}| x_i^\maxvar.
\end{equation*}
Dividing by $\tau$:
\begin{equation*}
  \frac{Z_j(\tau)}{\tau} \leq \frac{ V \prtr{ y_0^\maxvar - y_0^{(\epsilon)}} + E}{\epsilon \tau} + \frac{\sum_{i = 1}^n |a_{ji}| x_i^\maxvar}{\tau}.
\end{equation*}
When $V = 1/\varepsilon$ and $\tau = 1/\varepsilon^2$, it follows that
\begin{align*}
  \frac{Z_j(\tau)}{\tau} 
  & \leq \frac{ 1/\varepsilon \times \prtr{ y_0^\maxvar - y_0^{(\epsilon)}} + E}{\epsilon /\varepsilon^2} + \frac{\sum_{i = 1}^n |a_{ji}| x_i^\maxvar}{1/\varepsilon^2} \\
  & = O(\varepsilon).
\end{align*}
From \eqref{eq:rate_stable_bound}, it follows that
\begin{equation*}
  \sum_{i=1}^n a_{ji} \bar{x}_i(t) - b_j \leq O(\varepsilon).
\end{equation*}
Thus, under Assumption \ref{ass:staticslater}, using $O(1/\varepsilon^2)$ iterations ensures the time-averaged value of $\bar{y}_0(t)$ is within $O(\varepsilon)$ of the optimal value $y_0^\optvar$, and all constraints are within $O(\varepsilon)$ of being satisfied.  This is the $O(1/\varepsilon^2)$ tradeoff between computation and accuracy.


\subsection{Example}
For an example, we solved a small linear programming problem by using the max-weight and quadratic policies. The problem is
\begin{align*}
  \maximize \quad
  & 2 x_1 + ~x_2 \\
  \subjectto \quad
  & ~x_1 + ~x_2 \leq 4 \\
  & 5 x_1 + 3 x_2 \leq 15 \\
  & ~x_1 \quad\quad \leq 2.5 \\
  & 0 \leq x_1 \leq 10 \\
  & 0 \leq x_2 \leq 10.
\end{align*}

\begin{figure}
  \centering
  \includegraphics[scale=0.42]{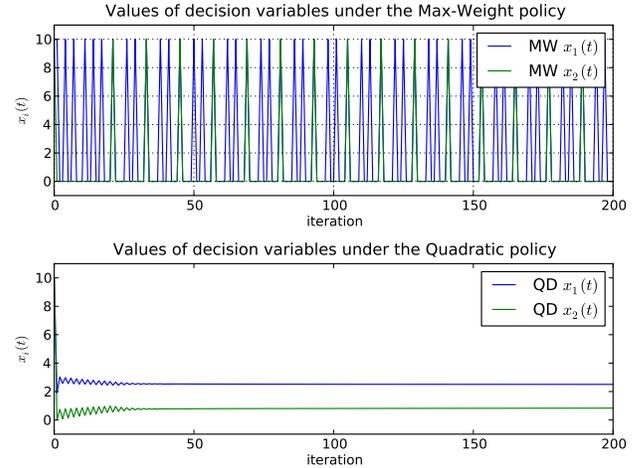}
  \caption{Comparison between max-weight and quadratic policies for solving linear programming}
  \label{fig:small_lp}
\end{figure}

The solution of this problem is $x_1 = 2.5, x_2 = 0.833$.  For both policies, the parameters are $V = 200$ and the number of iteration is $500$.  The values of decision variables $x_i(t)$ from both policies are shown in figure \ref{fig:small_lp}. The numerical values are show in table \ref{tlb:simple_problem}.
\begin{table}
  \caption{Numerical results from an example problem}
  \label{tlb:simple_problem}
  \centering
  \begin{tabular}{|c||c|c|c|}
    \hline\hline
          & Quadratic & Max-weight & Optimal \\
    \hline\hline
    $\bar{x}_1(500)$ & $2.531$ & $2.540$ & $2.500$ \\
    $\bar{x}_2(500)$ & $0.834$ & $0.820$ & $0.833$ \\
    \hline
    $x_1(500)$       & $2.500$ & $0.000$ & $2.500$ \\
    $x_2(500)$       & $0.833$ & $0.000$ & $0.833$ \\
    \hline\hline
  \end{tabular}
\end{table}

These time-averaged values of decision variables from both policies approach the optimal solution.  If number of iteration is increased, the precision is increased. Interestingly, the quadratic policy has a smooth property, as shown in Fig. \ref{fig:small_lp}, and that the intermediate decision values converge to an optimal solution before the time-averaged values does.

\section{Conclusion}
\label{sec:conclusion}
We studied information quality maximization in a system with uplink and single-hop relay capability which was done by designing queuing dynamic.  From Lyapunov optimization theory, we proposed a novel quadratic policy having a separable property, which leads to a distributed mechanism of format selection.  In comparison with the standard method, max-weight policy, our policy leads to an algorithm that reduces queue backlog by a significant constant.  This reduction also propagates and grows with the number of queues in the system.  We simulated the algorithm to verify correctness and behavior of the new policy.  In addition, we shows how the novel policy is applied to solve linear programs.

\bibliographystyle{IEEEtran}
\bibliography{Reference}

\end{document}

%% file: SuchaLib.tex
\newcommand{\set}[1]{ {\mathcal{#1}} }
\newcommand{\vect}[1]{ {\boldsymbol{#1}} }
\newcommand{\expect}[1]{ {\mathbb{E}\left\{ {#1} \right\}} }

\newcommand{\abs}[1]{ {\left| {#1} \right|} }

\newcommand{\maxvar}{ {(\text{max})} }

\newcommand{\actvar}{ {(\text{act})} }
\newcommand{\optvar}{ {(\text{opt})} }
\newcommand{\bestvar}{ {(\text{best})} }

\newcommand{\maximize}{ {\text{Maximize}} }
\newcommand{\minimize}{ {\text{Minimize}} }
\newcommand{\subjectto}{ {\text{Subject to}} }

\newcommand{\RealSet}{ {\mathbb{R}} }

\newcommand{\prts}[1]{ {\left[ {#1} \right]} }
\newcommand{\prtr}[1]{ {\left( {#1} \right)} }
\newcommand{\prtc}[1]{ {\left\{ {#1} \right\}} }

\newcommand{\prtbs}[1]{ {\biggl[ {#1} \biggr]} }

\newcommand{\control}{ {\alpha} }
\newcommand{\random}{ {\omega} }

\newcommand{\queue}{ {\boldsymbol{\Theta}} }
\newcommand{\defequiv}{\triangleq}

\DeclareMathOperator{\argmin}{argmin}
